\documentclass[12pt]{article}
\usepackage[amsmath]{e-jc_mod}
\usepackage{amssymb}
\usepackage{graphicx}
\usepackage{float}

\DeclareMathOperator{\lcm}{lcm}
\def\vec#1{\mathchoice{\mbox{\boldmath$\displaystyle\bf#1$}}
{\mbox{\boldmath$\textstyle\bf#1$}}
{\mbox{\boldmath$\scriptstyle\bf#1$}}
{\mbox{\boldmath$\scriptscriptstyle\bf#1$}}}

\def\a{\vec a}

\theoremstyle{plain}
\newtheorem{notation}[theorem]{Notation}
\dateline{Mar 31, 2022}{Aug 16, 2023}{TBD}
\Copyright{The authors. Released under the CC BY-ND license (International 4.0).}

\title{The generalized Frobenius problem via\\restricted partition functions}
\author{Kevin Woods\authornote{1}}

\authortext{1}{Department of Mathematics, Oberlin College, Oberlin, Ohio, USA (\email{Kevin.Woods@oberlin.edu})}

\begin{document}

\providecommand{\abs}[1]{\lvert#1\rvert}
\providecommand{\floor}[1]{\left\lfloor#1\right\rfloor}
\providecommand{\Z}{\mathbb{Z}}
\providecommand{\Zo}{\mathbb{Z}_{\ge 0}}
\providecommand{\R}{\mathbb{R}}
\providecommand{\N}{\mathbb{N}} \providecommand{\C}{\mathbb{C}}
\providecommand{\Q}{{\mathbb{Q}}} \providecommand{\x}{\mathbf{x}}
\providecommand{\y}{\mathbf{y}} \providecommand{\z}{\mathbf{z}}

\maketitle

\begin{abstract}
Given relatively prime positive integers, $a_1,\ldots,a_n$, the Frobenius number is the largest integer with no representations of the form $a_1x_1+\cdots+a_nx_n$ with nonnegative integers $x_i$. This classical value has recently been generalized: given a nonnegative integer $k$, what is the largest integer with at most $k$ such representations? Other classical values can be generalized too: for example, how many nonnegative integers are representable in at most $k$ ways? For sufficiently large $k$, we give formulas for these values by understanding the level sets of the restricted partition function (the function $f(t)$ giving the number of representations of $t$). Furthermore, we give the full asymptotics of all of these values, as well as reprove formulas for some special cases (such as the $n=2$ case and a certain extremal family from the literature). Finally, we obtain the first two leading terms of the restricted partition function as a so-called quasi-polynomial.
\end{abstract}

\medskip

\section{Introduction}
Given relatively prime positive integers, $a_1,\ldots,a_n$, we define the \emph{Frobenius number} to be the largest integer not contained in the semigroup
\[\left\{a_1x_1+\cdots+a_nx_n:\ x_i\in\Zo\right\}.\]
Formulas for some special cases have been known since at least Sylvester \cite{Sylvester} in the 1880's; for example, if $n=2$, the Frobenius number is $a_1a_2-a_1-a_2$.
See the Ram\'{\i}rez Alfons\'{\i}n text \cite{RA} for much more background.

More recently, Beck and Robins \cite{BR_frob} propose a generalization. While the classical Frobenius number is the largest integer that can be represented as a nonnegative integer combination of $a_1,\ldots,a_n$ in \emph{zero} ways,  we could instead take a fixed $k$ and look at integers that can be represented in exactly $k$ distinct ways. To be precise:

\begin{definition}
\label{def:g}
Given a vector $\a=(a_1,\ldots,a_n)$ of relatively prime positive integers and given $t\in\Zo$, define the \emph{restricted partition function}
\[f(\a; t)=\# (x_1,\ldots,x_n)\in\Zo^n:\ a_1x_1+\cdots+a_nx_n=t\]
to be the number of ways to represent $t$ by a nonnegative integer combination of the $a_i$. We write it as $f(t)$ when $\a $ is clear from context. Then define
\begin{itemize}
\item $g_{= k}$ to be the maximum $t\in\Zo$ such that $f(t)= k$ (the largest integer that can be represented in \emph{precisely} $k$ ways), if any such $t$ exist, and
\item $g_{\le k}$ to be the maximum $t\in\Zo$ such that $f(t)\le k$ (the largest integer that can be represented in \emph{at most} $k$ ways).
\end{itemize}
\end{definition}
\smallskip

The Frobenius number is $g_{=0}=g_{\le 0}$, but these numbers may differ for larger $k$:
\begin{example}(Shallit and Stankewicz \cite{SS}) For $\a=(8,9,15)$, we have $g_{=15}=169$, but $g_{\le 15}=g_{=14}=172$.
\end{example}
\smallskip

\begin{remark}
A consequence of Theorem \ref{thm:g} will be that $g_{=k}=g_{\le k}$, for all sufficiently large $k$.
\end{remark}
\smallskip

\begin{figure}[H]
\includegraphics[width=\textwidth]{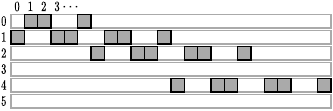}
\caption{The horizontal axis is $t=0,1,2,3,\ldots$ and the vertical axis is $f(t)$, in Example \ref{ex:interlace1}.}
\label{fig:interlace}
\end{figure}

\begin{example}
\label{ex:interlace1}
Take $\a=(3,4,6)$. Here is a table of $t$ and $f(t)$ for small $t$:
\footnotesize
\[\begin{array}{r|c|c|c|c|c|c|c|c|c|c|c|c|c|c|c|c|c|c|c|c|c|c|c|c|c} t& 0 & 1 & 2 & 3 & 4 & 5 & 6 & 7 & 8 & 9 & 10 & 11 & 12 & 13 & 14 & 15 & 16 & 17 & 18 & 19 & 20 & 21 & 22 & 23 & \cdots\\
\hline
f(t)  & 1 & 0 & 0 & 1 & 1 & 0 & 2 & 1 & 1 & 2 & 2 & 1 & 4 & 2 & 2 & 4 & 4 & 2 & 6 & 4 & 4 & 6 & 6 & 4 & \cdots
\end{array} 
\]
\normalsize
For example, $g_{=0}=5$ is the Frobenius number, and $g_{=2}=17$; the two representations of 17 are $17=3\cdot 1 + 4\cdot 2+6\cdot 1=3\cdot 3+4\cdot 2+6\cdot 0$. Except for $k=0$, which appears 3 times on this list of $f(t)$, values of $k$ seem to appear either 6 times ($k=1,2,4,\ldots$) or not at all $k=3,5,\ldots$. Figure \ref{fig:interlace} (inspired by Bardomero and Beck \cite[Figure 1]{BB}) illustrates how the level sets of $f(t)$ ``interlace'': the nonempty levels sets (except for $f(t)=0$) are translates of each other that eventually tile $\Z_{\ge 0}$.
\end{example}
\smallskip

In order to attack the generalized Frobenius problem, we will generalize Figure \ref{fig:interlace} and characterize how the level sets of $f(t)$ will interlace and how they will increase with $t$. We will make heavy use of the fact that $f(t)$ is a very ``nice'' function. In fact, it is a \emph{quasi-polynomial}:

\begin{definition}
A function $f:\Zo\rightarrow \Q$ is a \emph{quasi-polynomial of period $m$} if there exist polynomials $f_0,f_1,\ldots,f_{m-1}\in\Q[t]$ such that
\[f(t)=f_i(t),\text{ for }t\equiv i\bmod m.\]
The polynomials, $f_i$, are called the \emph{constituent polynomials} of $f$.
\end{definition}
\smallskip

The following folklore theorem shows that our $f$ is a quasi-polynomial:

\begin{proposition}
\label{prop:ehr}
Let $\a=(a_1,\ldots,a_n)$ be a vector of relatively prime positive integers. Then $f(\a;t)$ is a quasi-polynomial of period $m=\lcm(\a)=\lcm(a_1,\ldots,a_n)$. Furthermore, the leading term of all of the constituent polynomials is
\[\frac{1}{(n-1)!a_1\cdots a_n}t^{n-1}.\]
\end{proposition}
\smallskip

This proposition is apparently due to Issai Schur; see Wilf \cite[Section 3.15]{Wilf}, and we present a proof as part of Proposition \ref{prop:pf}.

Our first theorem will tell us exactly how to determine whether $f(s)=f(t)$, $f(s)>f(t)$, or $f(s)<f(t)$, for sufficiently large $s$ and $t$, and elucidate the structure of the output of $f$. First some notation:

\begin{notation}
For $1\le i\le n$, define $\a_{-i}=(a_1,\ldots,a_{i-1},a_{i+1},\ldots,a_n)$, so that, for example, $\gcd(\a_{-i})=\gcd(a_1,\ldots,a_{i-1},a_{i+1},\ldots,a_n)$.
\end{notation}

\smallskip

\begin{theorem}
\label{thm:f}
Let $\a=(a_1,\ldots,a_n)$ be a vector of relatively prime positive integers. For $1\le i\le n$, let $d_i=\gcd(\a_{-i})$, and let $p=d_1\cdots d_n$. Then
\begin{enumerate}
\item Let 
\[L=\left\{\sum_i a_ib_i:\ b_i\in \Z, 0\le b_i< d_i\right\}.\]
If $s\in\Zo$ and $\ell\in L$, then
\[f\left(sp+\ell\right)=f(sp).\]
(These will give the level sets of $f$, for sufficiently large $t$, all translates of $L$.)
\item Given $t\in\Zo$, there exists $s\in\Z$ and $\ell\in L$ such that
\[t=sp+\ell.\]
Furthermore, if $f(t)>0$, then $s\ge 0$.
(That is, Part (1) gives all of the level sets except for $f(t)=0$.)
\item For $1\le i\le n$, let
\[a'_i=\frac{a_i}{\prod_{j\ne i}d_i}\]
and $\a'=(a'_1,\ldots,a'_n)$. Then
\[f(\a; sp)=f(\a'; p),\]
for $s\in\Zo$. (This will be useful to simplify calculations of $f(sp)$, when $p>1$.)
\item For all sufficiently large $s\in\Zo$, 
\[f\big(\a;(s+1)p\big)>f(\a; sp).\]
(So these interlaced level sets will be broadly increasing with $t$.)
\end{enumerate}
\end{theorem}
\smallskip

\begin{example}
Continuing Example \ref{ex:interlace1} with $\a=(3,4,6)$, we can now better understand Figure \ref{fig:interlace}. Since $d_1=\gcd(4,6)=2$, $d_2=\gcd(3,6)=3$, and $d_3=\gcd(3,4)=1$, we have $p=2\cdot 3\cdot 1=6$. The set of values in $L$ are:
\begin{align*}3\cdot 0 + 4\cdot 0+6\cdot 0&=0,\quad 3\cdot 0 + 4\cdot 1+6\cdot 0=4,\quad 3\cdot 0 + 4\cdot 2+6\cdot 0=8,\\
3\cdot 1 + 4\cdot 0+6\cdot 0&=3,\quad 3\cdot 1 + 4\cdot 1+6\cdot 0=7,\quad 3\cdot 1 + 4\cdot 2+6\cdot 0=11.\end{align*}
Therefore, given $s\in\Zo$, $f(6s+\ell)$ will be identical for $\ell\in L=\{0,3,4,7,8,11\}$, which is exactly what we see in Figure \ref{fig:interlace}. Furthermore, the value of $f(6s)$ will eventually increase with $s$; in this example, it is increasing for all $s$: $f(0)=1$, $f(6)=2$, $f(12)=4$, $f(18)=6$, and so on. Except for $f(t)=0$, all values in the range $f(\Zo)$ will appear on this list (these interlaced translates of $L$ tile $\{t\in\Zo:\ f(t)> 0\}$).

Finally, $a'_1=3/3=1$, $a'_2=4/2=2$, and $a'_3=6/6=1$. One can check by hand that
\[f(\a; 6s)=f(\a'; s)=\begin{cases} \frac{s^2}{4}+s+1&\text{if $s$ is even,}\\  \frac{s^2}{4}+s+\frac{3}{4} &\text{if $s$ is odd.}\end{cases}\]
\end{example}
\smallskip

\begin{remark}
Since $f(t)$ is a quasi-polynomial of period $m=\lcm(\a)$ and we only need to look at values of $t$ that are multiples of $p$, we must compute $m/p$ of the constituent polynomials of $f$. In the above Example, $m/p=12/6=2$ and we need two polynomials.
\end{remark}
\smallskip

We now describe what this means for $g_{=k}$ and $g_{\le k}$, for sufficiently large $k$. We also describe some other quantities that often appear in both the classical and generalized Frobenius problem. Roughly, $g_{=k}$ and $g_{\le k}$ find the \emph{maximum} $t$ with a given property, but we might also want to find the \emph{minimum} such $t$, \emph{count} all such $t$, or even \emph{sum} all such $t$; \emph{generating functions} have also proven useful in studying these properties, so we analyze them too. The long list of precise definitions below --- and the parts of theorems pertaining to them --- can be skipped on first reading, in order to focus on $g_{=k}$ and $g_{\le k}$.

\begin{definition}
\label{def:frobs}
Let $\a=(a_1,\ldots,a_n)$ be a vector of relatively prime positive integers. For $k\in \Zo$, define
\begin{itemize}
\item $h_{= k}$ to be the minimum $t\in\Zo$ such that $f(t)= k$ (if any such $t$ exist),
\item $h_{\ge k}$ to be the minimum $t\in\Zo$ such that $f(t)\ge k$,
\item $c_{=k}$ to be the number of $t\in\Zo$ such that $f(t)= k$, 
\item $c_{\le k}$ to be the number of $t\in\Zo$ such that $f(t)\le k$,
\item $s_{=k}$ to be the sum of all $t\in\Zo$ such that $f(t)= k$,
\item $s_{\le k}$ to be the sum of all $t\in\Zo$ such that $f(t)\le k$,
\item $F_{= k}(x)$ to be the generating function
\[\sum_{t\in\Zo:\ f(t)= k}x^t,\]
\item $F_{\ge k}(x)$ to be the generating function
\[\sum_{t\in\Zo:\ f(t)\ge k}x^t.\]
\end{itemize}
\end{definition}
\smallskip

\begin{theorem}
\label{thm:g}
Let $\a=(a_1,\ldots,a_n)$ be a vector of relatively prime positive integers, and define $p,d_i$ as in Theorem \ref{thm:f}. Define $g_{=k}, g_{\le k}, h_{=k}, h_{\ge k}, c_{=k}, c_{\le k}, s_{=k},s_{\le k},$ $F_{= k}(x), F_{\ge k}(x)$ as in Definitions \ref{def:g} and \ref{def:frobs}. Then there are constants $C_1,C_2$ such that, for sufficiently large $s\in\Zo$,
\begin{align*}
g_{=f(sp)}=g_{\le f(sp)}&=sp+\sum_{i=1}^n(d_i-1)a_i,\\
h_{=f(sp)}=h_{\ge f(sp)}&=sp,\\
c_{=f(sp)}&=p,\\
c_{\le f(sp)}&=sp+C_1,\\
s_{=f(sp)}&=sp^2+\sum_{i=1}^n\frac{pa_i(d_i-1)}{2},\\
s_{\le f(sp)}&=\frac{1}{2}(sp)^2+\left(\frac{p+\sum_{i=1}^n a_i(d_i-1)}{2}\right)sp+C_2,\\
F_{= f(sp)}(x)&=x^{sp}\prod_i\frac{1-x^{d_ia_i}}{1-x^{a_i}},\\
F_{\ge f(sp)}(x)&=\frac{x^{sp}}{1-x^p}\prod_i\frac{1-x^{d_ia_i}}{1-x^{a_i}}.
\end{align*}
\end{theorem}
\smallskip

\begin{remark}
Let's call $k$ such that no $t$ has exactly $k$ representations ($c_{=k}=0$) \emph{trivial}. By Theorem \ref{thm:f},  the only \emph{nontrivial} $k$ are of the form $k=f(sp)$, and the values of such $g_{=f(sp)}$, etc., are given by the above theorem (for sufficiently large $s$). But this also gives us the values for (sufficiently large) \emph{trivial} $k$: for example,
$g_{\le k}=g_{\le f(sp)}$, if $f(sp)\le k < f\big((s+1)p\big)$.
\end{remark}
\smallskip

Notice that $g_{\le k}$ (like several of the other quantities) is of the form $g_{\le f(sp)}=sp+C$, where $C$ is a constant. That is, it is roughly the inverse of $f$. Writing $q_1(x)\sim q_2(x)$, if $\lim_{x\rightarrow\infty}q_1(x)/q_2(x)=1,$ Proposition \ref{prop:ehr} gives that
\[f(t)\sim\frac{1}{(n-1)!a_1\cdots a_n}t^{n-1}.\]
Therefore, if $k\sim f(sp)$ (in particular, if $f(sp)\le k< f\big((s+1)p\big)$), we have
\[sp\sim \big((n-1)!a_1\cdots a_nk\big)^{1/(n-1)},\]
and we immediately get the asymptotics of these functions of $k$:

\begin{corollary}
Given  a vector $\a=(a_1,\ldots,a_n)$ of relatively prime positive integers, let $p$ be the constant defined in Theorem \ref{thm:f}. Then (restricting to $k$ where the values are defined/nonzero)
\begin{itemize}
\item $g_{=k},g_{\le k},h_{=k},h_{\ge k},c_{\le k}$ \ $\sim$ \ $\big((n-1)!a_1\cdots a_nk\big)^{1/(n-1)}$,
\item $c_{=k}\sim p$,
\item $s_{=k}\sim p\big((n-1)!a_1\cdots a_nk\big)^{1/(n-1)}$,
\item $s_{\le k}\sim\frac{1}{2}\big((n-1)!a_1\cdots a_nk\big)^{2/(n-1)}$.
\end{itemize}
\end{corollary}
\smallskip

Fukshansky and Sch\"{u}rmann \cite{FS} give bounds for $g_{\le k}$, for sufficiently large $k$, matching these asymptotics, and Aliev, Fukshansky, and Henk \cite{AFH} find  bounds on $g_{\le k}$ that are good for all $k$. The asymptotics of the other quantities seem to be new here.

These quantities have already been calculated exactly for $n=2$, in Beck and Robins \cite{BR_frob} and Bardomero and Beck \cite{BB}. We will reproduce these results nicely using Theorem \ref{thm:f}:

\begin{proposition}
\label{prop:two}
Given relatively prime positive integers $a_1,a_2$, 
\begin{align*}
g_{= k}=g_{\le k}&=(k+1)a_1a_2-a_1-a_2,\\
\text{ for $k\ge 1$, }h_{= k}=h_{\ge k}&=(k-1)a_1a_2,\\
h_{=0}&=1 \text{ (unless $a_1=1$ or $a_2=1$)},\\
\text{ for $k\ge 1$, }c_{= k}&=a_1a_2,\\
c_{= 0}&= \frac{a_1a_2-a_1-a_2+1}{2},\\
c_{\le k}&=ka_1a_2+c_{=0},\\
\text{ for $k\ge 1$, }s_{= k}&=\frac{a_1a_2(2a_1a_2k-a_1-a_2)}{2},\\
s_{= 0}&=\frac{(a_1-1)(a_2-1)(2a_1a_2-a_1-a_2-1)}{12},\\
s_{\le k}&=\frac{a_1^2a_2^2}{2}k^2+\frac{a_1a_2(a_1a_2-a_1-a_2)}{2}k+s_{=0},\\
\text{ for $k\ge 1$, }F_{=k}(x)&= \frac{x^{(k-1)a_1a_2}\left(1-x^{a_1a_2}\right)^2}{\left(1-x^{a_1}\right)\left(1-x^{a_2}\right)},\\
F_{=0}(x)&=\frac{1}{1-x}-\frac{1-x^{a_1a_2}}{\left(1-x^{a_1}\right)\left(1-x^{a_2}\right)},\\
\text{ for $k\ge 1$, }F_{\ge k}(x)&=\frac{x^{(k-1)a_1a_2}\left(1-x^{a_1a_2}\right)}{\left(1-x^{a_1}\right)\left(1-x^{a_2}\right)},\\
F_{\ge 0}&=\frac{1}{1-x}.
\end{align*}

\end{proposition}
\smallskip
The formulas for $g_{=k},g_{\le k},h_{=k},h_{\ge k},c_{=k},c_{\le k}$ are due to (or immediately derivable from) \cite{BR_frob} and the formulas for $s_{=k},s_{\le k},F_{=k}(s),F_{\ge k}(x)$ are due to \cite{BB}. The $k=0$ cases were previously known: see Sylvester \cite{Sylvester} for $g_{=0}, c_{=0}$, Brown and Shiue \cite{Brown} for $s_{=0}$, and Sz\'{e}kely and Wormald \cite{SW} for $F_{=0}(x),F_{\ge 1}(x)$. Proposition \ref{prop:two} is an immediate corollary (the $n=2$ case) of Proposition \ref{prop:extremal} and Remark \ref{rem:0} below:

\begin{proposition}
\label{prop:extremal}
Let $d_1,\ldots,d_n$ be pairwise coprime positive integers, and let $a_i=\prod_{j\ne i}d_i$, for $1\le i\le n$. Let $p=d_1\cdots d_n$ and $\sigma=a_1+\cdots+a_n$. Other than $k=0$, the only nontrivial $k$ (that is, such that $c_{=k}> 0$) are $k=\binom{s+n-1}{n-1}$, for $s\in\Zo$, and we have
\begin{align*}
g_{= k}=g_{\le k}&=(s+n)p-\sigma,\\
h_{= k}=h_{\ge k}&= sp,\\
c_{= k}&=p,\\
c_{\le k}&=(s+1)p+\frac{(n-1)p-\sigma+1}{2},\\
s_{= k}&=\frac{p\big((2s+n)p-\sigma\big)}{2},\\
F_{= k}(x)&=\frac{x^{sp}\left(1-x^{p}\right)^{n}}{\left(1-x^{a_1}\right)\cdots \left(1-x^{a_n}\right)},\\
F_{\ge k}(x)&=\frac{x^{sp}\left(1-x^{p}\right)^{n-1}}{\left(1-x^{a_1}\right)\cdots \left(1-x^{a_n}\right)}.
\end{align*}

\end{proposition}
\smallskip

The formula for $g_{=k}=g_{\le k}$ was given in Beck and Kifer \cite{BK}. The other formulas seem to be new. If $n=2$, then $a_1=d_2$ and $a_2=d_1$ are generic relatively prime positive integers, and setting $k=\binom{s+1}{1}=s+1$ retrieves Proposition \ref{prop:two} for $k\ge 1$; the $k=0$ case is covered by the following remark:

\begin{remark}
\label{rem:0}
For $k=0$, Tripathi \cite{Tripathi} proved that
\[g_{=0}= (n-1)p-\sigma\quad\text{and}\quad c_{=0}=\frac{(n-1)p-\sigma+1}{2}.\]
These can be instead be obtained directly from $F_{\ge 1}(x)$ above, as follows: We have
\[F_{\ge 0}(x)=\sum_{t\in\Zo}x^t=\frac{1}{1-x}\quad\text{and}\quad F_{=0}(x)=F_{\ge 0}(x)-F_{\ge 1}(x).\]
Then $g_{=0}$ is the degree of $F_{=0}(x)$ as a polynomial and $c_{=0}=F_{=0}(1)$, which matches Tripathi's  \cite{Tripathi}  formulas.
One could compute $s_{=0}=F'_{=0}(1)$, which would also allow us to give a formula for $s_{\le k}$, but the answer seems a bit messy; however, $F'_{=0}(1)$ does match the $n=2$ value of $s_{=0}$ given in Proposition \ref{prop:two}.
\end{remark}
\smallskip

The following well-known lemma gives a useful recurrence and is worth highlighting here:
\begin{lemma}
\label{lem:recursion}
Given $t\in\Zo$, and given $i$ with $t\ge a_i$,
\[f(\a; t)=f(\a;t-a_i)+f(\a_{-i};t).\]
If we define $f(\a;t)=0$  for $t<0$ and $f(\emptyset;0)=1$, this recurrence holds for all $t\in\Z$.
\end{lemma}
The proof is immediate: the first term on the right-hand-side is the number of ways to represent $t$ with at least one $a_i$, and the second term is the number of ways to represent $t$ with no $a_i$'s.

Finally, we note that a partial fractions approach provides an alternative proof of Theorem \ref{thm:f}(4), and a standard proof of Proposition \ref{prop:ehr}. We include it here, in case it is useful. While the leading term of $f(\a; t)$ is well-known, this approach (together with Theorem \ref{thm:f}) also allows us to compute the second leading term(s) as well:

\begin{proposition}
\label{prop:pf}
Let $\a=(a_1,\ldots,a_n)$ be a vector of relatively prime positive integers, and let $m=\lcm(\a)$. For $1\le i\le n$, let $d_i=\gcd(\a_{-i})$, and let $p=d_1\cdots d_n$. Then
\begin{enumerate}
\item $f(\a; t)$ is a quasi-polynomial of period $m$, and the leading term of all of the constituent polynomials is
\[\frac{1}{(n-1)!a_1\cdots a_n}t^{n-1}.\]

\item If $d_i=1$ for all $i$, then the leading two terms of all of the constituent polynomials are
\[\frac{1}{(n-1)!a_1\cdots a_n}t^{n-1}+\frac{a_1+\cdots+a_n}{2(n-2)!a_1\cdots a_n}t^{n-2}.\]

\item For sufficiently large $s\in\Zo$, $f\big((s+1)p\big)>f(sp)$.
\end{enumerate}
\end{proposition}
\smallskip

\begin{remark}
Combining Proposition \ref{prop:pf}(2) and Theorem \ref{thm:f} allows us to compute the leading two terms even when $d_i>1$, though the second term will now depend on the constituent polynomial: given $t\in\Zo$, compute $r\in\Zo$ such that $t\equiv r\pmod p$ and $f(\a;t)=f(\a;t-r)$, using Theorem \ref{thm:f}(1) and (2) ($r$ depends only on $t$ mod $p$). Let $s\in\Z$ be such that $t=sp+r$, and then
\[f(\a;t)=f(\a;sp)=f(\a';s),\]
by Theorem \ref{thm:f}(3). The two leading terms of $f(\a';s)$ are given by Proposition \ref{prop:pf}(2), and then these can be used to compute the two leading terms of $f(\a;t)$ as a quasi-polynomial in $t$, by substituting $s=(t-r)/p$. The second leading term will depend on $t$ mod $p$.
\end{remark}
\smallskip

In the next section, we prove Theorem \ref{thm:f}, Theorem \ref{thm:g}, Proposition \ref{prop:extremal}, and Proposition \ref{prop:pf}. Then we conclude with some open questions.
\section{Proofs}

\begin{proof}[Proof of Theorem \ref{thm:f}]
Part 1 follows from the recurrence, Lemma \ref{lem:recursion}. In particular, we proceed by induction on $\ell=\sum_j b_j$. If all $b_j$ are zero, then this is trivially true: $f(sp+0)=f(sp)$. Now assume $b_i>0$, for some $i$. By Lemma \ref{lem:recursion} and the induction hypothesis,
\begin{align*}f\left(\a; sp+\sum_j a_jb_j\right)&=f\left(\a; sp+a_i(b_i-1)+\sum_{j\ne i}a_jb_j\right)+f\left(\a_{-i}; sp+\sum_j a_jb_j\right)\\
&=f\left(\a; sp\right)+f\left(\a_{-i}; sp+\sum_j a_jb_j\right).
\end{align*}
We need to show that $f\left(\a_{-i}; sp+\sum_j a_jb_j\right)=0$. Indeed, using the facts that $p$ and $a_j$ ($j\ne i$) are multiples of $d_i=\gcd(\a_{-i})$, that $a_i$ is relatively prime to $d_i$ (or else $\gcd(\a)>1$), and $b_i$ is not a multiple of $d_i$ (since $0<b_i<d_i$), we have
\[sp+\sum_j a_jb_j\equiv a_ib_i\not\equiv 0\pmod{d_i}.\]
Such a number cannot be represented as a combination of $\a_{-i}$, since $a_j$ ($j\ne i$) are multiples of $d_i$.
\\ \smallskip

Part 2 uses a standard number theory trick to compute $\ell=\sum_j b_j$. In particular, given $t\in\Zo$ let $b_i$ ($1\le i\le n$) be defined so that $0\le b_i<d_i$ and $b_i\equiv a_i^{-1}t\pmod{d_i}$ ($a_i$ is invertible mod $d_i$, since they are relatively prime). Since $a_j$ ($j\ne i$) is a multiple of $d_i$,
\[\sum_j a_jb_j\equiv a_ib_i\equiv t\pmod{d_i}.\]
Since $p=d_1\cdots d_n$ with the $d_i$  pairwise coprime (or else $\gcd(\a)>1$), the Chinese Remainder Theorem yields $\sum_j a_jb_j\equiv t\pmod p$. Let $s$ be the integer $\left(t-\sum_j a_jb_j\right)/p$, so that $t=sp+\sum_j a_jb_j$, as desired.

Now assume $f(t)>0$, and we need to prove $s\ge 0$. Recall that if we define $f(t)=0$ for $t<0$, then the recurrence in Lemma \ref{lem:recursion} applies for all $t\in\Z$, and therefore Part 1 (which only used that recurrence) holds for all $s\in\Z$. Then
\[f(sp)=f\left(sp+\sum_j a_jb_j\right)=f(t)>0,\]
which requires that $s\ge 0$, as desired.
\\ \smallskip

To prove Part 3, we must relate representations using $\vec a$ to representations using $\vec a'$. In particular, suppose $sp=\sum_j a_jx_j$ ($x_j\in \Zo$) is a representation of $sp$ by $\a$. For each $i$, $p$ and $a_j$ ($j\ne i$) are multiples of $d_i$, and so
\[a_ix_i\equiv \sum_j a_jx_j=sp\equiv 0\pmod{d_i}.\] Since $a_i$ and $d_i$ are relatively prime, $x_i$ must be a multiple of $d_i$. Let $y_i\in\Zo$ be such that $x_i=d_iy_i$. Then
\[sp=\sum_i a_ix_i=\sum_i \left(\prod_{j\ne i}d_i\right)a'_i\cdot d_iy_i=p\sum_i a'_iy_i,\]
So $s=\sum_i a'_iy_i$ is a representation of $s$ by $\a'$. Conversely, given any representation $s=\sum_i a'_iy_i$  ($y_i\in\Zo$) by $\a'$, $sp=\sum_i a_i(d_iy_i)$ is a representation of $sp$ by $\a$. Therefore $f(\a; sp)=f(\a'; p)$, as desired.
\\ \smallskip

Part 4 requires a deeper understanding of the function $f(t)$.  First, we assume without loss of generality (by Part 3) that $d_i=1$ for all $i$, so we are trying to prove that $f(s+1)>f(s)$, for sufficiently large $s\in\Zo$. The complication is that $f(s)$ and $f(s+1)$ are evaluated on different constituent polynomials of $f$, and it seems like these might ``jump around.''  We use the recurrence, Lemma \ref{lem:recursion}, to show that $f(s)$ and $f(s+1)$ can both be related to the same $f(s-q)$ and therefore to each other, and this relation will entail that $f(s+1)-f(s)$ is eventually positive.

Indeed, we know that all sufficiently large integers can be represented by $\a$. In particular, let $q\in\Zo$ be such that $q$ and $q+1$ are both representable; that is, $q=\sum_i a_ix_i$ and $q+1=\sum_i a_iy_i$ for $x_i,y_i\in\Zo$. Take $s\in\Zo$ sufficiently large (in particular, take $s\ge q$). We will use Lemma \ref{lem:recursion} repeatedly to relate both $f(s)$ and $f(s+1)$ to $f(s-q)$. Let's start by applying the recursion $x_1$ times on $f(\a; s)$, using $i=1$:
\begin{align*}
f(\a; s)&=f(\a; s-a_1)+f(\a_{-1};s)=\\
	&=f(\a; s-2a_1)+f(\a_{-1};s-a_1)+f(\a_{-1};s)=\cdots\\
&=f(\a;s-a_1x_1)+\sum_{j=0}^{x_1-1}f(\a_{-1};s-ja_1).
\end{align*}
Now apply the recursion $x_2$ times with $i=2$, and so on, and we get constants (independent of $s$) $u_{ij}\in\Zo$ such that
\[ f(\a; s)=f(\a;s-q)+\sum_{i=1}^n\sum_{j=0}^{x_i-1}f(\a_{-i};s-u_{ij}).\]
Now if we do the same thing for $f(\a; s+1)$, applying the recursion $y_1$ times with $i=1$ and so forth, we get constants $w_{ij}\in\Zo$ such that
\[ f(\a; s+1)=f\big(\a;s+1-(q+1)\big)+\sum_{i=1}^n\sum_{j=0}^{y_i-1}f(\a_{-i};s+1-w_{ij}).\]
Subtracting the two equations, the term $f(\a;s-q)=f\big(\a;s+1-(q+1)\big)$ cancels, and we are left with
\[f(\a; s+1) - f(\a; s) = \sum_{i=1}^n\sum_{j=0}^{y_i-1}f(\a_{-i};s+1-w_{ij}) - \sum_{i=1}^n\sum_{j=0}^{x_i-1}f(\a_{-i};s-u_{ij}),\]
and we want to show that this quantity is (eventually) positive. By Proposition \ref{prop:ehr}, $f(\a_{-i};s)$ is a quasi-polynomial with leading term
\[\frac{1}{(n-2)!\prod_{j\ne i}a_j}s^{n-2}\]
(note that we are using that $\gcd(\vec a_{-i})=d_i=1$).
Therefore $f(\a; s+1) - f(\a; s)$ is a quasi-polynomial with leading coefficient (on $s^{n-2}$)
\begin{align*}
\sum_{i=1}^n&\sum_{j=0}^{y_i-1}\frac{1}{(n-2)!\prod_{j\ne i}a_j} - \sum_{i=1}^n\sum_{j=0}^{x_i-1}\frac{1}{(n-2)!\prod_{j\ne i}a_j}\\
&=\sum_{i=1}^n \frac{a_iy_i}{(n-2)!a_1\cdots a_n}-\sum_{i=1}^n \frac{a_ix_i}{(n-2)!a_1\cdots a_n}\\
&=\frac{1}{(n-2)!a_1\cdots a_n}\left[\sum_{i=1}^na_iy_i-\sum_{i=1}^na_ix_i\right]\\
&=\frac{1}{(n-2)!a_1\cdots a_n}\big((q+1)-q\big)\\
&=\frac{1}{(n-2)!a_1\cdots a_n}.
\end{align*}
Since this is a positive leading term, $f(\a; s+1) - f(\a; s)$ will eventually be positive, as desired.

\end{proof}
\smallskip

\begin{proof}[Proof of Theorem \ref{thm:g}]
Let $s$ be sufficiently large. By Theorem \ref{thm:f}, the set of $t$ with $f(t)=f(sp)$ is exactly
\[\left\{sp+\sum_i a_ib_i:\ 0\le b_i<d_i\right\}.\]
We simply need to check what that means for all of our different values:\\

\smallskip
\noindent The largest element of this set occurs at $b_i=d_i-1$ for all $i$, so
\[g_{=f(sp)}=sp+\sum_{i=1}^n(d_i-1)a_i.\]
The smallest element of this set occurs at $b_i=0$ for all $i$, so
\[h_{=f(sp)}=sp.\]
The number of elements in this set is
\[c_{=f(sp)}=d_1\cdots d_n=p.\]
The sum of the elements in this set is
\begin{align*}s_{=f(sp)}&=\sum_{b_1=0}^{d_1-1}\cdots\sum_{b_n=0}^{d_n-1}\left(sp+\sum_{i=1}^n a_ib_i\right)\\
&=d_1\cdots d_n sp + \sum_{i=1}^{n}\left(\prod_{j\ne i}d_j\cdot\sum_{b_i=0}^{d_i-1}a_ib_i \right)\\
&=sp^2+ \sum_{i=1}^{n}\left(\prod_{j\ne i}d_j\cdot\frac{a_id_i(d_i-1)}{2}\right)\\
&=sp^2+\sum_{i=1}^n\frac{pa_i(d_i-1)}{2}.
\end{align*}
The generating function for this set is
\begin{align*}
F_{= f(sp)}(x)&=\sum_{b_1=0}^{d_1-1}\cdots\sum_{b_n=0}^{d_n-1}x^{sp+\sum_{i=1}^n a_ib_i}\\
&=x^{sp}\prod_i \left(1+x^{a_i}+\cdots+x^{(d_i-1)a_i}\right)\\
&=x^{sp}\prod_i\frac{1-x^{d_ia_i}}{1-x^{a_i}}.
\end{align*}

Since $f(sp)$ is an increasing function of $s$ (for sufficiently large $s$), we have that $g_{\le f(sp)}=g_{=f(sp)}$ and $h_{\ge f(sp)}=h_{=f(sp)}$. To compute $c_{\le f(sp)}$, we have to worry about small $k$. In particular, $c_{=0}$ might not be $p$ (see Example \ref{ex:interlace1}), and it is possible that $f(rp)=f(r'p)$ for distinct (small) $r,r'$ so that some $c_{=k}$ is a nontrivial multiple of $p$. But we will have (for sufficiently large $s$) that
\[c_{\le f(sp)}=c_{=0}+\sum_{r=0}^s p=sp+C_1,\]
where $C_1$ is a constant. Similarly, we may compute
\begin{align*}
s_{\le f(sp)}&=s_{=0}+\sum_{r=0}^s \left(rp^2+\sum_{i=1}^n\frac{pa_i(d_i-1)}{2}\right)\\
&=s_{=0}+\frac{s^2p^2+sp^2}{2}+(s+1)\sum_{i=1}^n\frac{pa_i(d_i-1)}{2}\\
&=\frac{1}{2}(sp)^2+\left(\frac{p+\sum_{i=1}^n a_i(d_i-1)}{2}\right)sp+C_2,
\end{align*}
where $C_2$ is a constant.
Finally,
\begin{align*}
F_{\ge f(sp)}(x)&=\sum_{r=s}^{\infty}F_{=f(rp)}(x)\\
&=\sum_{r=s}^{\infty}x^{rp}\prod_i\frac{1-x^{d_ia_i}}{1-x^{a_i}}\\
&=\frac{x^{sp}}{1-x^p}\prod_i\frac{1-x^{d_ia_i}}{1-x^{a_i}}.
\end{align*}
\end{proof}
\smallskip

\begin{proof}[Proof of Proposition \ref{prop:extremal}]
In this setting, Theorem \ref{thm:g}(3) drastically simplifies the calculation of $f$, allowing us to make quick work of the rest. In particular,  note that the $d_i$ as defined in the proposition are indeed $d_i=\gcd(\a_{-i})$, as required to apply Theorem \ref{thm:g}. By Theorem \ref{thm:f}(1) and (2), we may concentrate on $f(sp)$ for $s\in\Zo$. So let $s\in \Zo$ be given, and let $k=\binom{s+n-1}{n-1}$. We have (in the notation of Theorem \ref{thm:f}(3))
\[a'_i=\frac{a_i}{\prod_{j\ne i}d_i}=1,\]
for all $i$. Then by Theorem \ref{thm:f}(3),
\[f(\a;sp)=f(\a';s)=f\big((1,\ldots,1); s\big)=\binom{s+n-1}{n-1}=k\]
(this calculation is a classical combinatorics problem on compositions: $f\big((1,\ldots,1); s\big)$ is the number of ways to write $s=x_1+\cdots+x_n$, where $x_i\in\Zo$, which is the number of ways to shuffle $s$ identical ``stars'' and $n-1$ identical ``bars'').
Now we may simply apply Theorem \ref{thm:g}, and using that $p=a_id_i$:

\begin{align*}
g_{=f(sp)}=g_{\le f(sp)}&=sp+\sum_{i=1}^n(d_i-1)a_i\\
&=sp+\sum_{i=1}^n \left(p-a_i\right)\\
&=(s+n)p-\sigma,\\
h_{=f(sp)}=h_{\ge f(sp)}&=sp,\\
c_{=f(sp)}&=p,\\
s_{=f(sp)}&=sp^2+\sum_{i=1}^n\frac{pa_i(d_i-1)}{2}\\
&=sp^2+\sum_{i=1}^n\frac{p^2-pa_i}{2}\\
&=\frac{2sp^2}{2}+\frac{np^2-p\sigma}{2}\\
&=\frac{p\big((2s+n)p-\sigma\big)}{2},\\
F_{= f(sp)}(x)&=x^{sp}\prod_i\frac{1-x^{d_ia_i}}{1-x^{a_i}}\\
&=\frac{x^{sp}\left(1-x^{p}\right)^{n}}{\left(1-x^{a_1}\right)\cdots \left(1-x^{a_n}\right)},\\
F_{\ge f(sp)}(x)&=\frac{x^{sp}}{1-x^p}\prod_i\frac{1-x^{d_ia_i}}{1-x^{a_i}}\\
&=\frac{x^{sp}\left(1-x^{p}\right)^{n}}{\left(1-x^p\right)\left(1-x^{a_1}\right)\cdots \left(1-x^{a_n}\right)}\\
&=\frac{x^{sp}\left(1-x^{p}\right)^{n-1}}{\left(1-x^{a_1}\right)\cdots \left(1-x^{a_n}\right)}.
\end{align*}
Finally, using that $c_{=0}=\frac{(n-1)p-\sigma+1}{2}$ from Tripathi \cite{Tripathi},
\begin{align*}c_{\le k}&=c_{=0}+\sum_{r=0}^s c_{=\binom{r+n-1}{n-1}}\\
&=\frac{(n-1)p-\sigma+1}{2}+\sum_{r=0}^s p\\
&=(s+1)p+\frac{(n-1)p-\sigma+1}{2}.\\
\end{align*}
\end{proof}
\smallskip

\begin{proof}[Proof of Proposition \ref{prop:pf}]
We will be brief, since much of this is classical; see Wilf's text \cite[Section 3.15]{Wilf}, for example. Define $G(x)=\sum_{t=0}^{\infty} f(\a;t)x^t.$ We see that
\[G(x)=(1+x^{a_1}+x^{2a_1}+\cdots)\cdots(1+x^{a_n}+x^{2a_n}+\cdots)=\frac{1}{\prod_i (1-x^{a_i})}.\]
We will use the partial fraction expansion of $G(x)$ to get our results. All of the poles of $G$ are $m$th roots of unity, where $m=\lcm(\a)$. One pole is $x=1$, of order $n$. Label the other roots of unity by $\zeta_j$, for $1\le j<m$, and suppose $\zeta_j$ is a pole of order $b_j$. Then the partial fraction expansion of $G(x)$ yields that there exist $C_\ell,D_{j\ell}\in \Q$ such that
\[G(x)=\sum_{\ell=1}^{n}\frac{C_\ell}{(1-x)^\ell}+\sum_{j=1}^{m-1}\sum_{\ell=1}^{b_j}\frac{D_{j\ell}}{(1-x/\zeta_j)^\ell}.\]

Suppose $\zeta_j$ is a primitive $r$th root of unity. Then a term $\frac{D_{j\ell}}{(1-x/\zeta_j)^\ell}$, if expanded out as a product of geometric series,  contributes a degree $\ell-1$ quasi-polynomial of period $r$ to $f(t)$. Summed together, we will have a period $m$ quasi-polynomial. We can see that $\zeta_j$ is a root of exactly those $1-x^{a_i}$ such that $r$ divides $a_i$; therefore, it will be a pole of order $b_j=\abs{i:\ r\text{ divides }a_i}$.

Since $\gcd(\a)=1$, we must have $b_j\le n-1$, and so the only degree $n-1$ piece will come from 
\[C_n/(1-x)^n=\sum_{t=0}^{\infty}C_n\binom{t+n-1}{n-1}x^t\]
(the $t$th coefficient in the power series will be the number of ways to write $t=c_1+\cdots+c_n$ with $c_i\in\Zo$, the same classic combinatorics problem as in the proof of Proposition \ref{prop:extremal}).

Furthermore, if $d_i=1$ for all $i$, then no $r>1$ can divide $n-1$ of the $a_i$, and so $b_j\le n-2$, and the only degree $n-1$ and $n-2$ pieces will come from
\[C_n/(1-x)^n+C_{n-1}/(1-x)^{n-1}=\sum_{t=0}^{\infty}\left(C_n\binom{t+n-1}{n-1}+C_{n-1}\binom{t+n-2}{n-2}\right)x^t.\]

Noting that
\[C_n=(1-x)^{n}G(x)\Big|_{x=1}\quad\text{and}\quad C_{n-1}=\frac{d}{dx}(1-x)^{n}G(x)\Big|_{x=1},\]
we compute that
\[C_n=\frac{1}{a_1\cdots a_n}\quad\text{and}\quad C_{n-1}=\frac{a_1+\cdots+a_n-n}{2a_1\cdots a_n},\]
and we can compute that
\begin{align*}C_n&\binom{t+n-1}{n-1}+C_{n-1}\binom{t+n-2}{n-2}\\
&=\frac{1}{(n-1)!a_1\cdots a_n}t^{n-1}+\frac{a_1+\cdots+a_n}{2(n-2)!a_1\cdots a_n}t^{n-2}+\text{ lower order terms}.
\end{align*}
This gives the first leading term of $f(t)$, in general, and the first two leading terms when $d_i=1$ for all $i$, and so Parts (1) and (2) are proved.
\smallskip

To prove Part (3), Theorem\ref{thm:g}(3) allows us to assume without loss of generality that $d_i=1$ for all $i$, and we want to prove that $f(s+1)>f(s)$ for sufficiently large $s$. Indeed, the leading term of $f(s+1)-f(s)$, when expanded out as a quasi-polynomial using Part (2), is
\[\frac{1}{(n-2)!a_1\cdots a_n}s^{n-2}.\]
Since this is a positive leading term, $f(s+1)-f(s)$ must eventually be positive, as desired.

\end{proof}

\section{Open Questions}
\begin{question}
We have made no effort to quantify what \emph{sufficiently large} means in any of these theorems, but probably one can, since $f(t)$ is so ``well-behaved'' here. What bounds can we give for when the results  hold?
\end{question}
\smallskip

\begin{question}
The $n=2$ case is well understood (see Proposition \ref{prop:two}), and finding formulas for $n\ge 4$ seems very difficult even in the $k=0$ case. It seems possible that there are interesting formulas when $n=3$, however. For example, when $n=3$ and  $k=0$, there are reasonable formulas (see Ram\'{\i}rez Alfons\'{\i}n  \cite[Chapter 2]{RA}, and, for a generating function approach, see Denham \cite{Denham}). Are there interesting formulas for $n=3$ and general $k$?
\end{question}
\smallskip

\begin{question}
Let $P\subseteq\R^n$ be a $d$-dimensional polytope whose vertices are rational, and let $m$ be the smallest integer such that the vertices of $mP$ ($P$ dilated by a factor of $m$) are integers. Then Ehrhart \cite{Ehrhart} proves that $f(t)=\abs{tP\cap\Z^n}$ is a quasi-polynomial of period $m$ (see the Beck and Robins text \cite{BR_book} for many more details). This is a generalization of our problem, as taking $P$ to be the convex hull of $\vec e_i/a_i$ ($1\le i\le n$), where $\vec e_i$ is $i$th standard basis vector, yields the Frobenius $f(t)$. One can define $g_{\le k}$, and so forth, using this new $f$, and Aliev, De Loera, and Louveaux \cite{ADL}  study structural and algorithmic results related to this. Do some of the results of this current paper generalize to that more general setting?
\end{question}
\smallskip

\begin{question}
What can we say about the computational complexity of computing $g_{\le k}$, $c_{\le k}$, and so forth? If $n$ is not fixed, then Ram\'{\i}rez Alfons\'{\i}n  \cite{RA_NP} shows that even computing $g_{= 0}$ is NP-hard. On the other hand, if $n$ is fixed, then Kannan \cite{Kannan} shows that  $g_{=0}$ can be computed in polynomial time, and Barvinok and Woods \cite{BW} show that  $c_{=0}$ and other quantities can be computed in polynomial time. Generalizing, Aliev, De Loera, and Louveaux \cite{ADL} show that, for fixed $n$ and $k$, $g_{\le k}$ and other quantities can be computed in polynomial time, even in the general setting of  $f(t)=\abs{tP\cap\Z^n}$.

This leaves the open question: Can these quantities be computed in polynomial time if $n$ is fixed, but $a_1,\ldots,a_n$ and $k$ are the input? Nguyen and Pak \cite{NP_Ehrhart} prove that this is NP-hard in the general setting of $f(t)=\abs{tP\cap\Z^n}$, disproving a conjecture from \cite{ADL}. However, to do this, they construct a polytope $P\subseteq\R^6$ whose $f(t)$ varies wildly across the constituent polynomials, which is not true for our Frobenius $f(t)$ (see Theorem~\ref{thm:f}).

When $k$ is sufficiently large, Theorem \ref{thm:g} applies: For any given $t$, we can compute $f(t)$ in polynomial time, using the result of Barvinok \cite{Barvinok} that $\abs{P\cap\Z^n}$ can be computed in polynomial time for fixed $n$; then binary search allows us to find $s$ such that $f(sp)\le k < f\big((s+1)p\big)$, and Theorem \ref{thm:g} gives us $g_{=k}$. But what if $k$ is bigger than a constant but not ``sufficiently large'' for Theorem \ref{thm:g} to hold?

\end{question}


\begin{thebibliography}{10}

\bibitem{ADL}
Iskander Aliev, Jes\'{u}s~A. De~Loera, and Quentin Louveaux.
\newblock Parametric polyhedra with at least {$k$} lattice points: their
  semigroup structure and the {$k$}-{F}robenius problem.
\newblock In {\em Recent trends in combinatorics}, volume 159 of {\em IMA Vol.
  Math. Appl.}, pages 753--778. Springer, 2016.

\bibitem{AFH}
Iskander Aliev, Lenny Fukshansky, and Martin Henk.
\newblock Generalized {F}robenius numbers: bounds and average behavior.
\newblock {\em Acta Arith.}, 155(1):53--62, 2012.

\bibitem{BB}
Leonardo Bardomero and Matthias Beck.
\newblock Frobenius coin-exchange generating functions.
\newblock {\em Amer. Math. Monthly}, 127(4):308--315, 2020.

\bibitem{BW}
Alexander Barvinok and Kevin Woods.
\newblock Short rational generating functions for lattice point problems.
\newblock {\em J. Amer. Math. Soc.}, 16(4):957--979, 2003.

\bibitem{Barvinok}
Alexander~I. Barvinok.
\newblock A polynomial time algorithm for counting integral points in polyhedra
  when the dimension is fixed.
\newblock {\em Math. Oper. Res.}, 19(4):769--779, 1994.

\bibitem{BK}
Matthias Beck and Curtis Kifer.
\newblock An extreme family of generalized {F}robenius numbers.
\newblock {\em Integers}, 11:A24, 6, 2011.

\bibitem{BR_frob}
Matthias Beck and Sinai Robins.
\newblock A formula related to the {F}robenius problem in two dimensions.
\newblock In {\em Number theory ({N}ew {Y}ork Seminar 2003)}, pages 17--23.
  Springer, New York, 2004.

\bibitem{BR_book}
Matthias Beck and Sinai Robins.
\newblock {\em Computing the continuous discretely}.
\newblock Undergraduate Texts in Mathematics. Springer, New York, 2007.

\bibitem{Brown}
Tom~C. Brown and Peter Jau-Shyong Shiue.
\newblock A remark related to the {F}robenius problem.
\newblock {\em Fibonacci Quart.}, 31(1):32--36, 1993.

\bibitem{Denham}
Graham Denham.
\newblock Short generating functions for some semigroup algebras.
\newblock {\em Electron. J. Combin.}, 10:R36, 2003.

\bibitem{Ehrhart}
Eugene Ehrhart.
\newblock Sur les polyedres homothetiques bordes a n dimensions.
\newblock {\em Comptes Rendus Hebdomadaires des seances de l'academie des
  sciences}, 254(6):988, 1962.

\bibitem{FS}
Lenny Fukshansky and Achill Sch\"{u}rmann.
\newblock Bounds on generalized {F}robenius numbers.
\newblock {\em European J. Combin.}, 32(3):361--368, 2011.

\bibitem{Kannan}
Ravi Kannan.
\newblock Lattice translates of a polytope and the {F}robenius problem.
\newblock {\em Combinatorica}, 12(2):161--177, 1992.

\bibitem{NP_Ehrhart}
Danny Nguyen and Igor Pak.
\newblock On the number of integer points in translated and expanded polyhedra.
\newblock {\em Discrete Comput. Geom.}, 65(2):405--424, 2021.

\bibitem{RA_NP}
Jorge~L. Ram\'{\i}rez~Alfons\'{\i}n.
\newblock Complexity of the {F}robenius problem.
\newblock {\em Combinatorica}, 16(1):143--147, 1996.

\bibitem{RA}
Jorge~L. Ram\'{\i}rez~Alfons\'{\i}n.
\newblock {\em The {D}iophantine {F}robenius problem}, volume~30 of {\em Oxford
  Lecture Series in Mathematics and its Applications}.
\newblock Oxford University Press, Oxford, 2005.

\bibitem{SS}
Jeffrey Shallit and James Stankewicz.
\newblock Unbounded discrepancy in {F}robenius numbers.
\newblock {\em Integers}, 11:A2, 8, 2011.

\bibitem{Sylvester}
James~J Sylvester.
\newblock Mathematical questions with their solutions.
\newblock {\em Educational times}, 41(21):171--178, 1884.

\bibitem{SW}
L\'{a}zl\'{o}~A. Sz\'{e}kely and Nicholas~C. Wormald.
\newblock Generating functions for the {F}robenius problem with {$2$} and {$3$}
  generators.
\newblock {\em Math. Chronicle}, 15:49--57, 1986.

\bibitem{Tripathi}
Amitabha Tripathi.
\newblock On a linear {D}iophantine problem of {F}robenius.
\newblock {\em Integers}, 6:A14, 6, 2006.

\bibitem{Wilf}
Herbert~S. Wilf.
\newblock {\em generatingfunctionology}.
\newblock Academic Press, Inc., Boston, MA, second edition, 1994.

\end{thebibliography}
\end{document}